\documentclass[a4paper,11pt]{amsart}

\usepackage[margin=1in,includehead,includefoot]{geometry}
\usepackage[OT4]{fontenc}
\usepackage{amsthm,amsmath,amssymb}
\usepackage{enumitem}
\usepackage[pdftitle={Triangle-free geometric intersection graphs with no large independent sets},
            pdfauthor={Bartosz Walczak},
            colorlinks=true,linkcolor=black,citecolor=black,filecolor=black,urlcolor=black]{hyperref}

\newtheorem*{theorem}{Theorem}

\renewenvironment{enumerate}{\begin{enumorig}[label=\textup{(\roman*)}, noitemsep, topsep=1.5mm plus 1.5mm, leftmargin=*, widest=iii]}{\end{enumorig}}

\def\calP{\mathcal{P}}
\def\calI{\mathcal{I}}
\def\setR{\mathbb{R}}

\let\leq\leqslant
\let\geq\geqslant
\let\setminus\smallsetminus

\def\size#1{\lvert #1\rvert}

\makeatletter
\let\old@setaddresses\@setaddresses
\def\@setaddresses{\bigskip\bgroup\parindent 0pt\let\scshape\relax\old@setaddresses\egroup}
\makeatother

\title{Triangle-free geometric intersection graphs with~no~large~independent sets}

\author{Bartosz Walczak}

\address{Theoretical Computer Science Department, Faculty of Mathematics and Computer Science, Jagiellonian University, Krak\'ow, Poland}
\email{walczak@tcs.uj.edu.pl}

\thanks{A journal version of this paper appeared in \emph{Discrete Comput.\ Geom.}, 53(1):221--225, 2015.}
\thanks{The author was supported by Polish National Science Center grant 2011/03/B/ST6/01367 and Swiss National Science Foundation grant 200020-144531.}

\begin{document}

\begin{abstract}
It is proved that there are triangle-free intersection graphs of line segments in the plane with arbitrarily small ratio between the maximum size of an independent set and the total number of vertices.
\end{abstract}

\maketitle

\section{Introduction}

Pawlik et al.\ \cite{PKK+14} proved that there are triangle-free intersection graphs of line segments in the plane with arbitrarily large chromatic number.
The graphs they construct have independent sets containing more than $1/3$ of all the vertices.
It has been left open whether there is a constant $c>0$ such that every triangle-free intersection graph of $n$ segments in the plane has an independent set of size at least $cn$.
Fox and Pach \cite{FoP12} conjectured a much more general statement, that $K_k$-free intersection graphs of curves in the plane have linear-size independent sets, for every~$k$.
This would imply a well-known conjecture that $k$-quasi-planar graphs (graphs drawn in the plane so that no $k$ edges cross each other) have linearly many edges \cite{PSS96}, which is proved up to $k=4$ \cite{Ack09}.

In this note, I resolve the independent set problem in the negative, proving the following strengthening of the result of Pawlik et al.

\begin{theorem}
There are triangle-free segment intersection graphs with arbitrarily small ratio between the maximum size of an independent set and the total number of vertices.
\end{theorem}

The constructions presented in the next two sections give rise to triangle-free intersection graphs of $n$ segments in the plane with maximum independent set size $\Theta(n/\log\log n)$.

\section{Construction}

Pawlik et al.\ \cite{PKK+14} construct, for $k\geq 1$, a triangle-free graph $G_k$ and a family $\calP_k$ of subsets of $V(G_k)$, called \emph{probes}, with the following properties:
\begin{enumerate}
\item\label{item:probe-no} $\size{\calP_k}=\smash{2^{2^{k-1}-1}}$,
\item\label{item:probe-indep} every member of $\calP_k$ is an independent set of $G_k$,
\item\label{item:coloring} for every proper coloring of the vertices of $G_k$, there is a probe $P\in\calP_k$ such that at least $k$ colors are used on the vertices in $P$.
\end{enumerate}
They are built by induction on $k$, as follows.
The graph $G_1$ has just one vertex $v$, and $\calP_1$ has just one probe $\{v\}$.
For $k\geq 2$, first, a copy $(G,\calP)$ of $(G_{k-1},\calP_{k-1})$ is taken.
Then, for every probe $P\in\calP$, another copy $(G_P,\calP_P)$ of $(G_{k-1},\calP_{k-1})$ is taken.
There are no edges between vertices from different copies.
Finally, for every probe $P\in\calP$ and every probe $Q\in\calP_P$, a new vertex $d_Q$ connected to all vertices in $Q$, called the \emph{diagonal} of $Q$, is added.
The resulting graph is $G_k$.
The family of probes $\calP_k$ is defined by
\begin{equation*}
\calP_k=\bigl\{P\cup Q\colon P\in\calP\text{ and }Q\in\calP_P\bigr\}\cup\bigl\{P\cup\{d_Q\}\colon P\in\calP\text{ and }Q\in\calP_P\bigr\}.
\end{equation*}
It is easy to check that the graph $G_k$ is indeed triangle-free and the conditions \ref{item:probe-no}--\ref{item:coloring} are satisfied for $(G_k,\calP_k)$---see \cite{PKK+14} for details.
It is also shown in \cite{PKK+14} how the graph $G_k$ is represented as a segment intersection graph.

I will show that there is an assignment $w_k$ of positive integer weights to the vertices of $G_k$ with the following properties:
\begin{enumerate}
\setcounter{enumi}{3}
\item\label{item:weight} the total weight of $G_k$ is $\frac{k+1}{2}\cdot\smash{2^{2^{k-1}-1}}$,
\item\label{item:indep} for every independent set $I$ of $G_k$, the number of probes $P\in\calP_k$ such that $P\cap I\neq\emptyset$ is at least the weight of $I$.
\end{enumerate}
Once this is achieved, the proof of the theorem of this paper follows easily.
Namely, it follows from \ref{item:probe-no} and \ref{item:indep} that every independent set $I$ of $G_k$ has weight at most $\smash{2^{2^{k-1}-1}}$.
We can take the representation of $G_k$ as a segment intersection graph and replace every segment representing a vertex $v\in V(G_k)$ by $w_k(v)$ parallel segments lying very close to each other, so as to keep the property that any two segments representing vertices $u,v\in V(G_k)$ intersect if and only if $uv\in E(G_k)$.
It follows from \ref{item:weight} that the family of segments obtained this way has size $\frac{k+1}{2}\cdot\smash{2^{2^{k-1}-1}}$, while every independent set of its intersection graph has size at most $\smash{2^{2^{k-1}-1}}$.

The assignment $w_k$ of weights to the vertices of $G_k$ is defined by induction on $k$, following the inductive construction of $(G_k,\calP_k)$.
The weight of the only vertex of $G_1$ is set to $1$.
This clearly satisfies \ref{item:weight} and \ref{item:indep}.
For $k\geq 2$, let $G$, $\calP$, $G_P$, $\calP_P$ and $d_Q$ be defined as in the inductive step of the construction of $(G_k,\calP_k)$.
Let $p=|\calP_{k-1}|=\smash{2^{2^{k-2}-1}}$.
The weights $w_k$ of the vertices of $G$ are their original weights $w_{k-1}$ in $G_{k-1}$ multiplied by $p$.
The weights $w_k$ of the vertices of every $G_P$ are equal to their original weights $w_{k-1}$ in $G_{k-1}$.
The weight $w_k$ of every diagonal $d_Q$ is set to $1$.
It remains to prove that \ref{item:weight} and \ref{item:indep} are satisfied for $(G_k,\calP_k,w_k)$ assuming that they hold for $(G_{k-1},\calP_{k-1},w_{k-1})$.

The proof of \ref{item:weight} is straightforward:
\begin{equation*}
w_k(G_k)=w_k(G)+\sum_{P\in\calP}\bigl(w_k(G_P)+\size{\calP_P}\bigr)=2pw_{k-1}(G_{k-1})+p^2=\tfrac{k+1}{2}\cdot 2^{2^{k-1}-1}.
\end{equation*}

For the proof of \ref{item:indep}, let $I$ be an independent set in $G_k$.
Let $\calI=\{P\in\calP\colon P\cap I\neq\emptyset\}$.
For every probe $P\in\calP$, define
\begin{align*}
\calI_P&=\{Q\in\calP_P\colon Q\cap I\neq\emptyset\},&
\calP'_P&=\{P\cup Q\colon Q\in\calP_P\}\cup\{P\cup\{d_Q\}\colon Q\in\calP_P\},\\
D_P&=\{d_Q\colon Q\in\calP_P\},&
\calI'_P&=\{P'\in\calP'_P\colon P'\cap I\neq\emptyset\}.
\end{align*}
By the induction hypothesis, we have
\begin{equation*}
w_k(V(G)\cap I)\leq p\size{\calI},\qquad\qquad w_k(V(G_P)\cap I)\leq\size{\calI_P}.
\end{equation*}

Suppose $P\in\calI$.
It follows that $(P\cup Q)\cap I\neq\emptyset$ and $(P\cup\{d_Q\})\cap I\neq\emptyset$ for every $Q\in\calP_P$.
Hence $\size{\calI'_P}=\size{\calP'_P}=2p$.
Moreover, we have $d_Q\notin I$ whenever $Q\in\calI_P$, because $d_Q$ is connected to all vertices in $Q$, one of which belongs to $I$.
Hence
\begin{equation*}
w_k(V(G_P)\cap I)+w_k(D_P\cap I)\leq\size{\calI_P}+\size{\calP_P\setminus\calI_P}=\size{\calP_P}=p.
\end{equation*}

Now, suppose $P\in\calP\setminus\calI$.
If $Q\in\calI_P$, then $(P\cup Q)\cap I\neq\emptyset$, $d_Q\notin I$ (by the same argument as above), and $(P\cup\{d_Q\})\cap I=\emptyset$.
If $Q\in\calP_P\setminus\calI_P$, then $(P\cup Q)\cap I=\emptyset$, and $(P\cup\{d_Q\})\cap I\neq\emptyset$ if and only if $d_Q\in I$.
Hence
\begin{equation*}
w_k(V(G_P)\cap I)+w_k(D_P\cap I)\leq\size{\calI_P}+\size{D_P\cap I}=\size{\calI'_P}.
\end{equation*}

To conclude, we gather all the inequalities and obtain
\begin{equation*}
\begin{split}
w_k(I)&=w_k(V(G)\cap I)+\sum_{P\in\calP}\bigl(w_k(V(G_P)\cap I)+w_k(D_P\cap I)\bigr)\\
&\leq p\size{\calI}+\sum_{P\in\calI}p+\sum_{P\in\calP\setminus\calI}\size{\calI'_P}=\sum_{P\in\calI}\size{\calI'_P}+\sum_{P\in\calP\setminus\calI}\size{\calI'_P}=\sum_{P\in\calP}\size{\calI'_P}.\end{split}
\end{equation*}

\section{Improved construction}

Pawlik et al.\ \cite{PKK+14} define also a graph $\tilde G_k$, which arises from $(G_k,\calP_k)$ by adding, for every probe $P\in\calP_k$, a diagonal $d_P$ connected to all vertices in $P$.
This is the smallest triangle-free segment intersection graph known to have chromatic number greater than $k$.
Define the assignment $\tilde w_k$ of weights to the vertices of $\tilde G_k$ so that $\tilde w_k$ is equal to $w_k$ on the vertices of $G_k$ and $\tilde w_k(d_P)=1$ for every $P\in\calP_k$.
Let $I$ be an independent set in $\tilde G_k$.
Let $\calI=\{P\in\calP_k\colon P\cap I\neq\emptyset\}$.
Hence $d_P\notin I$ for $P\in\calI$.
It follows that
\begin{gather*}
\tilde w_k(I)=w_k(V(G_k)\cap I)+\size{\{d_P\colon P\in\calP_k\}\cap I}\leq\size{\calI}+\size{\calP_k\setminus\calI}=\size{\calP_k}=2^{2^{k-1}-1},\\
\tilde w_k(\tilde G_k)=w_k(G_k)+\size{\calP_k}=\tfrac{k+3}{2}\cdot 2^{2^{k-1}-1}.
\end{gather*}
The graph $\tilde G_k$ is the smallest one for which I can prove that it has a weight assignment such that the ratio between the maximum weight of an independent set and the total weight is at most $\smash{\frac{2}{k+3}}$.
It is not difficult to prove (e.g.\ using weak LP duality) that the assignment of weights $\tilde w_k$ to the vertices of $\tilde G_k$ is optimal (gives the least ratio) for this particular graph.

Both constructions give rise to triangle-free intersection graphs of $n$ segments in the plane with maximum independent set size $\Theta(n/\log\log n)$.
On the other hand, it follows from the result of McGuinness \cite{McG00} that every triangle-free intersection graph of $n$ segments has chromatic number $O(\log n)$ and maximum independent set size $\Omega(n/\log n)$.

\section{Other geometric shapes}

It is known that the graphs $G_k$ and $\tilde G_k$ have intersection models by many other geometric shapes, for example, L-shapes, axis-parallel ellipses, circles, axis-parallel square boundaries \cite{PKK+13} or axis-parallel boxes in $\setR^3$ \cite{Bur65}.
The result of this paper can be extended to those models for which every geometric object $X$ representing a vertex of the intersection graph can be replaced by many pairwise disjoint objects approximating $X$. 
This is possible, for example, for intersection graphs of L-shapes, circles or axis-parallel square boundaries, but not for intersection graphs of axis-parallel ellipses or axis-parallel boxes in $\setR^3$.
The problem whether triangle-free intersection graphs of the latter kind of shapes have linear-size independent sets remains open.

\section*{Acknowledgment}

I thank Michael Hoffmann for helpful discussions.

\bibliographystyle{plain}
\bibliography{independent}

\end{document}